\documentclass[12pt]{article}
\usepackage{amsmath,amsfonts,amsthm,amscd}
\usepackage [mathscr]{eucal}

\begin{document} 
\parskip=6pt
\title{Hochschild Cohomology and Twisted Complexes on Complex Manifolds}
\author{Yue Lin L. Tong \thanks{Research partially supported by Mathematics 
Division of NCTS at Taipei.}
\and I-Hsun Tsai \thanks{Research partially supported by a grant from NSC 
of Taiwan.}}\date{}

\maketitle

\begin{abstract}
We use the theory of twisted resolutions and twisted complexes to give a proof of Kontsevich's assertion that
Yoneda product and cup product are preserved in a canonical isomophism
$$
{Ext}^.(X\times X;{\cal O}_\Delta,{\cal O}_\Delta)\cong \oplus H^.(X,\wedge^.\,T_X)$$ where $X$ is a complex manifold and $\Delta$ is the diagonal in $X\times X$.
\end{abstract}

\begin{center}
{\bf{Introduction}}\\
\end{center}

\par The categorical mirror symmetry $[K1]$ has been an impetus for the considerable 
interest recently in the Hochschild cohomology of algebraic varieties
 $[C,K2,S,W1,Y]$. One of the themes is to unify the different possible 
definitions of
Hochschild cohomology, whether one takes it to be the $Ext$ groups, 
$Ext^.(X\!\times\! X; 
{\cal O}_\Delta,{\cal O}_\Delta)$, where $\Delta$ is diagonal in $X\times X$, or 
hypercohomology of
the standard Hochschild cochain complex, or variations thereof. Then using the 
degeneration of
spectral sequence for $Ext$, one gets via the Hochschild-Kostant-Rosenberg 
isomorphism a
decomposition:

\begin{equation}
Ext^.(X\!\times\! X;{\cal O}_\Delta,{\cal O}_\Delta)\cong \oplus 
H^.(X,{\Lambda}^.T_X)
\tag{0.1}
\end{equation}

\par An interesting problem, inspired by Kontsevich's assertion at 
the end of
$[K2]$, is the existence of a natural isomorphism between the two sides of (0.1) 
which
preserves the Yoneda product in $Ext^.$ and the cup product in 
$H^.(X,{\Lambda}^.T_X)$. In a
recent paper $[C]$, it is conjectured that such an isomorphism must involve a 
correction of
the $HKR$ isomorphism by the square root of the Todd class in the right hand 
side of (0.1).

\par The purpose of this paper is to show that the theory of twisted resolution 
and twisted
complexes $[TT1,2]$ is ideally suited to unify the different versions of 
Hochschild
cohomology on any complex manifold, and furthermore we show that $HKR$ 
isomorphism already
preserves products in (0.1), so no correction is needed. In comparison to the 
statement of
Kontsevich $[K2]$, whose approach presumably involves his deformation 
quantization techniques,
no such techniques are used here.

\par The sheaf ${\cal O}_\Delta$in $X\times X$ has concrete local Koszul 
resolutions in any
coordinate neighborhood of $\Delta$ in $X\times X$. Of course such local 
resolutions do not
patch up to a global resolution of ${\cal O}_\Delta$ in $X\!\times\! X$. In 
$[TT1,2]$, a technique is
introduced which builds a twisted resolution of ${\cal O}_\Delta$ from such 
local resolutions.
Using this data one can define twisted complexes to represent various cohomology 
functors on
$X\!\times\! X$ with coefficients in ${\cal O}_\Delta$, in particular the 
Hochschild cohomology
groups. The reliance on local resolutions overcomes a key difficulty which 
arises when the bar
resolution is used to represent $Ext^.(X\!\times\! X;{\cal O}_\Delta,{\cal 
O}_\Delta)$. The bar
resolution is only defined on affine open sets near the diagonal in $X\times X$. 
For example
the differentials are generally not defined at $(z,\zeta)\in X\!\times\! X$ if 
$z$ and $\zeta$ are not in the same affine open
set.

\par The restriction of local Koszul complex and its dual to $\Delta$ readily 
yield the sheaves $\Omega^._X$ and $\wedge^.T_X$, so the twisted complexes 
provide a particularly direct link between the two sides of (0.1). In fact the 
theory developed in $[TT2]$ already suffices to give an isomorphism in (0.1) 
which is
compatible with the products on both sides. However some work is required to 
show that the
isomorphism coincides with the $HKR$ isomorphism. For this purpose we
construct in \S3 a chain map between local bar and Koszul resolutions. The 
standard $HKR$ map
emerges natually when the chain map is restricted to $\Delta$. In \S 1 \&2 we 
recall the theory of twisted resolutions
and twisted complexes and make some adjustments for present applications. We 
have tried to
make the exposition reasonably self contained. In \S4 we represent Hochschild homology by the cohomology of a 
twisted
complex which naturally leads to $Tor.(X\!\times\! X;{\cal O}_\Delta,{\cal 
O}_\Delta)$. The details
here are mostly paralled to that of Hochschild cohomology developed in \S1-3.

\par Twisted resolutions also provide a convenient framework to consider deformations of coherent sheaves. We
hope to discuss related matters elsewhere.

\par We would like to thank the Mathematics Division of NCTS at Taipei for 
support on this
research.

\bigskip

\noindent{\bf \S1. Twisted Resolutions and Twisted Complexes}
\bigskip

We recall the construction of twisted resolutions and twisted complexes for the 
sheaf ${\cal
O}_\Delta[TT1,2]$. Let $X$ be a complex manifold of dimension $n$, and ${\cal 
U}=\{U_\alpha\}$ be a locally
finite Stein open cover of $X$ for which there are complex analytic coordinates 
$z^i_\alpha$ on each
$U_\alpha$. Let $p_1$, $p_2$ be projections $X\!\times\! X\to X$ onto the first 
and second
factors. On $U_\alpha \times U_\alpha$ write the local coordinates 
$p^*_1z^i_\alpha$,
$p^*_2z^i_\alpha$ as $z^i_\alpha$ and $\zeta^i_\alpha$ respectively. Take ${\cal 
W}$ as the open
cover of $X\times X$ consisting of the cover ${\cal U}\times {\cal U}=\{U_\alpha 
\times
U_\alpha;U_\alpha \in {\cal U}\}$ together with a Stein open cover ${\cal 
V}=\{V_\mu \}$ of the
complement of the diagonal.

\par Given an open set $W_\mu \in {\cal W}$ let $F_\mu$ be the free ${\cal 
O}_{W_\mu}$ module of rank
$n$, $F_\mu=\sum_{i=1}^n {\cal O}_{W_\mu}e^i$ and $K_\mu^{-q}=\Lambda^q(F_\mu)$, 
$0\le q\le n$. Also
let ${\check F}_\mu=\sum_{i=1}^n {\cal O}_{W_\mu}{\check e}^i$ be the dual 
module of $F_\mu$, and ${\check
K}_\mu^q=\Lambda ^q({\check F}_\mu)$, $0\le q\le n$. \newline 
If $W_\alpha=U_\alpha \times U_\alpha \in {\cal U}\times {\cal U}$ 
we have a differential

\begin{equation}
(d_K)_\alpha:K_\alpha^{-q} \to K_\alpha^{-q+1} 
\tag{1.1}
\end{equation}

\noindent given by contraction by $(z^1_\alpha-\zeta^1_\alpha)\check 
e^1+..+(z^n_\alpha
-\zeta^n_\alpha)\check e^n$, while a differential $(d_{\check K})_\alpha:{\check 
K}_\alpha^q \to {\check
K}_\alpha^{q+1}$ is wedge product with the same element. If $W_\mu \in {\cal 
V},(d_K)_\mu$ is
contraction by ${\check e}_1$ while $(d_{\check K})_\mu$ is wedge product by 
${\check e}_1$.
\bigskip

\par The Koszul complexes $K^.$ give local resolutions of the sheaf of the 
diagonal ${\cal O}_\Delta$
in $X\!\times\! X$. In fact chain homotopies on $U_\alpha \times 
U_\alpha:P_\alpha:K_\alpha^{-q} \to
K_\alpha^{-q-1}$ are given by

\begin{equation}
P_\alpha(fe^I)=\sum_{j=1}^n\left(\int_0^1 t^{|I|}\,\frac{\partial f}{\partial
z_\alpha^j}(\zeta_\alpha +t(z_\alpha-\zeta_\alpha),\zeta_\alpha)dt\right)e^j 
\wedge e^I  
\tag{1.2}
\end{equation}

\noindent where $I=(i_1..,i_q)$ is an increasing multi-index, and $|I|$ its 
length. We have $[TT1,
(9,19)]$.
\begin{equation}
(d_K)_\alpha P_\alpha + P_\alpha(d_K)_\alpha = 1- res  
\tag{1.3}
\end{equation}

\centerline{$P_\alpha^2=0$}
\noindent where $res:K_\alpha^{-q} \to K_\alpha^{-q}$ is zero except when $q=0$, 
and in that case for $f
\in K_\alpha^0$, 
$$
res(f)(z,\zeta) = f(z,z)
$$
\noindent If $W_\mu \in {\cal V}$, $P_\mu:K_\mu^{-q} \to K_\mu^{-q-1}$ is wedge 
product by $e^1$, then

\begin{equation}
(d_K)_\mu P_\mu + P_\mu(d_K)_\mu = 1  
\tag{1.4}
\end{equation}

\par The local resolutions$K^._\alpha \to {\cal O}_{(U_\alpha \times U_\alpha) 
\cap \Delta}$ cannot be
expected to patch up to a global resolution of ${\cal O}_\Delta$ in a 
neighborhood of $\Delta$ in $X
\times X$. Nevertheless it is possible to build a global twisted resolution from 
which one can represent
the various cohomological functors on $X\! \times \! X$ with coefficients in 
${\cal O}_\Delta$. \newline
\noindent In $W_{\alpha_0}\cap...\cap W_{\alpha_p} = W_{\alpha_0..\alpha_p}$ \ \ 
 we set

\begin{equation}
Hom^q(K^._{\alpha_p},K^._{\alpha_0}) = \oplus_i 
Hom(K^{-i}_{\alpha_p},K^{-i+q}_{\alpha_0})\cong \oplus_i
K^{-i+q}_{\alpha_0} \otimes {\check K}^i_{\alpha_p}  
\tag{1.5}
\end{equation}

\noindent The differential acting on $f \in K_{\alpha_0}^{-i+q} \otimes {\check 
K}^i_{\alpha_p}$ is
\begin{equation}
\begin{aligned}
(d_{Hom})_{\alpha_0..\alpha_p}(f)&= (d_K)_{\alpha_0}f 
+(-1)^{q+1}f(d_K)_{\alpha_p}\\
&=(d_K)_{\alpha_0}f +
(-1)^{i+q}(d_{\check K})_{\alpha_p} f 
\end{aligned}  
\tag{1.6}
\end{equation}

\noindent so that
\begin{equation}
d_{Hom}(fg)=(d_{Hom} f)g + (-1)^{|f|} f(d_{Hom} g) 
\tag{1.7}
\end{equation}

\noindent where $|f|$ is degree of $f$ in $Hom^.$. A homotopy operator $P_H$ is 
defined on $Hom^.$ by
using the tensor product representation in (1.5):
\begin{equation}
(P_H)_{\alpha_o..\alpha_p}(f) = \sum_{i \ge 0}(-1)^i P_{\alpha_o} 
\left(\sigma(d_{\check K})_{\alpha_p} P_{\alpha
_o}\right)^i(f) 
\tag{1.8}
\end{equation}

\noindent where $\sigma$ is the sign appearing in (1.6). Then by $[TT1, (8,17), 
TT2, (1,10)]$
\begin{equation}
d_{Hom} P_H + P_H d_{Hom} = 1 - r  
\tag{1.9}
\end{equation}

\noindent where $$r=\sum_{i\ge 0}(-1)^i(P_H \sigma d_{\check K})^i \, res$$.

\par Let $C^p({\cal W}, K^q)$ denote the group of ${\check C}\!ech \, p$ 
cochains whose values $c_{\alpha_
0..\alpha_p}$ lie in $K^q_{\alpha_0}|_{W_{\alpha_0..\alpha_p}}$. Similarly 
$C^p({\cal W}, Hom^r(K^.,
K^.))$ denote ${\check C}ech \, p$ cochains taking value in 
$Hom^r(K^._{\alpha_p}, K^._{\alpha_0})|_{W_{\alpha
_0..\alpha_p}}$. An associative cup product on $C^.({\cal W}, Hom^.(K^., K^.))$ 
as well as an action of
$C^.({\cal W}, Hom^.(K^.,K^.))$ on $C^.({\cal W}, K^.)$ are defined by 

\begin{equation}
(f^{pq} \cdot g^{rs})_{\alpha_0..\alpha_{p+r}}= (-1)^{qr} 
f^{pq}_{\alpha_0..\alpha_p}g^{rs}_{\alpha_p..\alpha_
{p+r}}  
\tag{1.10}
\end{equation}

\par The ${\check C}\!ech$ coboundary operators $\delta$ are given by the 
formulas
\begin{equation}
(\delta f)_{\alpha_0..\alpha_{p+1}}= \sum^p_{i=1}(-1)^i f_{\alpha_0..{\hat
\alpha}_i..\alpha_{p+1}}\bigm|_{ W_{\alpha_0..\alpha_{p+1}}}  
\tag{1.11}
\end{equation}

\noindent for $f \in C^p({\cal W}, {Hom}^.(K^., K^.))$ and 
$$(\delta c)_{\alpha_0...\alpha_{p+1}}=\sum^{p+1}_{i=1}(-1)^i c_{\alpha_0..{\hat
\alpha_i}..\alpha_{p+1}} \bigm|_{ W_{\alpha_0..\alpha_{p+1}}}
$$

\noindent for $c \in C^p({\cal W},K^.)$. It follows readily that
\begin{equation}
\delta (f \cdot g)= (\delta f)\cdot g+(-1)^{deg\, f}f \cdot (\delta g)  
\tag{1.12}
\end{equation}

\noindent where $\cdot$ is the product (1.10) and $deg f$ refers to the total 
degree $({\check C \!
ech}+{Hom}^.)$. These operators satisfy $\delta ^2 =0$, but they are different 
from standard ${\check C \!
ech}$ coboundary operators where the sum in (1.11) is over $0 \le i \le {p+1}$. 
To include the missing terms
we need to generalize the notion of transition functions of a complex of vector 
bundles. We consider the
following chain maps between the local Koszul complexes. Since $K^0_\alpha={\cal 
O}_{W_\alpha}$ define
$A^0_{\alpha \beta}:K^0_\beta \to K^0_\alpha$ to be just the identity map. The 
$A^0_{\alpha \beta}$ extends
to a chain map of the complexes by the formula.
\begin{equation}
A^._{\alpha \beta}=A^0_{\alpha \beta} - (P_H)_{\alpha \beta}(d_{{\check 
K}_\beta}A^0_{\alpha \beta}) 
\tag{1.13}
\end{equation}

\noindent Clearly $A^q_{\alpha \beta}=A^0_{\alpha \beta}$ when $q=0$, and 
$A^._{\alpha \beta}$ gives a chain
map, i.e. $A^._{\alpha \beta}$ is a $d_{Hom}$ cocycle from the equation, using 
(1.8):
$$
d_{Hom}\, A^._{\alpha \beta}= d_{{\check K}_ \beta} \, A^0_{\alpha \beta}- 
(1-P_H d_{Hom}-r)d_{{\check K}_
\beta} \, A^0_{\alpha \beta}=0, 
$$
\noindent since $res(d_{{\check K}_ \beta} \, A^0_{\alpha \beta})=0$.

\par Let $a^{1,0} \in C^1({\cal W}, {Hom}^0(K^., K^.))$ be given by 
$a^{1,0}_{\alpha \beta}=A^._{\alpha
\beta}$. Then since these are not true transition functions $a^{1,0}_{\alpha 
\beta}\, a^{1,0}_{\beta
\gamma}- a^{1,0}_{\alpha \gamma}$ does not vanish. However this difference is 
chain homotopic to zero, and thus gives
rise to 
$$
a^{i,-i+1} \in C^i({\cal W}, {Hom}^{-i+1}(K^., K^.))
$$
$0 \le i \le n$, where $a^{0,1}=d_{Hom}$, $a^{1,0}$ as above, and inductively,
\begin{equation}
a^{i+1,-i}=(-1)^i P_H (\delta a^{i,-i+1}+\sum_{\scriptstyle{r+ s=i-1 \atop r,s 
\ge 0}}a^{r+1,-r}\cdot
a^{s+1,-s}) 
\tag{1.14}
\end{equation}

\noindent Let $a=\sum_{i \ge 0}a^{i,-i+1}$, then $a$ has total degree 1 and it 
satisfies the twisting cochain
equation:  $[TT2,\S 2][OTT, \S 1]$ .
\begin{equation}
\delta a + a \cdot a = 0  
\tag{1.15}
\end{equation}

\noindent The data $({\cal W},K^., a)$ is called a twisted resolution of ${\cal 
O}_\Delta$ on $X\! \times\! X$.
This generalizes the 1 cocycle condition for transition functions of a global 
complex of vector bundles. The
equations (1.12) and (1.15) enable one to define differential $D_a$ on 
$C^.({\cal W},K^.)$ by
\begin{equation}
D_a \, c=\delta c+a \cdot c
\tag{1.16}
\end{equation}

\noindent and similarly differential $D_{a,a}$ on $C^.({\cal W}, 
{Hom}^.(K^.,K^.))$ \newline
by
\begin{equation}
D_{a,a}\,f= \delta f + a\cdot f + (-1)^{deg\, f+1}f\cdot a  
\tag{1.17}
\end{equation}

\noindent We denote the singly graded complexes by $C^._a({\cal W},K^.)$ and 
$C^._{a,a}({\cal
W},{Hom}^.(K^.,K^.))$ respectively. Thus 
$$
C^p_a({\cal W},K^.)=\bigoplus_{i\ge 0} C^{p+i}({\cal W}, K^{-i})
$$
\noindent with the differential $D_a$ and similarly for $C^p_{a,a}({\cal W}, 
{Hom}^.(K^.,K^.))$. These
differentials are compatible with products and actions (1.10).

\bigskip

\noindent{\bf \S2. Cohomology of Twisted Complexes}
\bigskip

Consider the restriction of the twisting cochain $a$ to the diagonal. From (1.1) 
it follows that
$a^{0,1}|_\Delta =0$. from $[TT2, p.53 (9),(6.1)]$ it follows that for $i \ge 2$
\begin{equation}
a^{i,-i+1}\bigm| _\Delta=0  
\tag{2.1}
\end{equation}
\noindent For the remaining component
\begin{equation}
a^{1,0}_{\alpha \beta}\bigm| _\Delta =\bigwedge^.\left(\frac{\partial 
z^i_{\beta}}{\partial z^j_{\alpha}}\right) 
\tag{2.2}
\end{equation}

\noindent This follows from the statement $[TT2,(3.10)]$, we fill in some 
details here as an exercise in the
formulas in \S1. From (1.8), (1.13) it follows that
\begin{equation}
A^{-1}_{\alpha \beta}= -P_\alpha(d_{{\check K}_\beta}\,A^0_{\alpha 
\beta})=P_\alpha(A^0_{\alpha
\beta}\,d_{{\check K}_\beta})
\tag{2.3}
\end{equation}

\noindent thus $A^{-1}_{\alpha \beta}(e^i)=P_\alpha(z^i_\beta-\zeta^i_\beta)$, 
and from (1.2)
$$
A^{-1}_{\alpha \beta}(e^i)\bigm| _\Delta=\sum_j\frac{\partial 
z^i_{\beta}}{\partial z^j_{\alpha}}e^j
$$

\noindent Now assume $A^{-k}_{\alpha \beta}|_\Delta$ acts as 
$\bigwedge^k\left(\frac{\partial
z^i_{\beta}}{\partial z^j_{\alpha}}\right)$ for $k<q$ and consider a muti index 
$I=(i_1,..,i_q)$. We denote by
$I_r$ the muti index $(i_1,..,{\hat i}_r,..,i_q)$.

\begin{equation}
\begin{aligned}
 A^{-q}_{\alpha \beta}(e^I)&=P_\alpha A^{-(q-1)}_{\alpha 
\beta}(d_{K_\beta}e^I)\\
&=P_\alpha \,A^{-(q-1)}_{\alpha
\beta}\left(\sum^q_{r=1}(-1)^{r-1}(z^{i_r}_\beta-\zeta^{i_r}_\beta)e^{I_r}\right
)\\
&=P_\alpha 
\left(\sum^q_{r=1}(-1)^{r-1}(z^{i_r}_\beta-\zeta^{i_r}_\beta)A^{-(q-1)}_{\alpha
\beta}(e^{I_r})\right) 
\end{aligned}
\tag{2.4}
\end{equation}


\noindent It is clear that in $A^{-q}_{\alpha \beta}(e^I)|\Delta$, the nonzero 
terms come from $P_\alpha$
differentiating the $(z^{i_r}_\beta-\zeta^{i_r}_\beta)$. Thus when restricted to 
$\Delta$ a summand in (2.4) is
given by 
$$
{(-1)^{r-1}\over{|I|}}\, \sum_j {{\partial z^{i_r}_{\beta}} \over \partial 
z^k_{\alpha}} 
\, e^k \bigwedge
\bigwedge^{q-1}({\partial z^i_{\beta}\over \partial 
z^j_{\alpha}})(e^{I_r})={1\over q}\, \bigwedge ^q({\partial z^i_{\beta}\over 
\partial z^j_{\alpha}})(e^I)$$

\noindent where we have used induction hypothesis and formula (1.2). This 
finishes the proof of (2.2). 
\par So when
restricted to $\Delta$, Koszul complexes may be identified with the sheaf of 
holomorphic forms.
\begin{equation}
K^{-q}_\alpha \bigm |_ \Delta \cong \Omega^q_{U_\alpha}  
\tag{2.5}
\end{equation}

Similarly the dual Koszul complex restricts to sheaf of multi tangent fields.
\begin{equation}
{\check K}^q_\alpha \bigm |_ \Delta \cong \bigwedge ^q \, T_{U_\alpha}  
\tag{2.6}
\end{equation}

\par The twisted complexes $C^._a$ and $C^._{a,a}$ have natual filtrations given 
by ${\check C}\! ech$
degrees which are preserved by differentials $D_a$ and $D_{a,a}$ respectively.

\bigskip
\noindent{{\bf{\underbar{Proposition}} (2.7)}}\bigskip
\begin{eqnarray*}
H^.({Hom}^.(K^.,K^.)) &\cong &\underbar {\it {Ext}}^._{{\cal O}_{X \times 
X}}({\cal O}_\Delta,{\cal O}_\Delta)\\
H^.(C^._{a,a}({\cal W},{Hom}^.(K^.,K^.))) &\cong& {Ext}^.(X \times X;{\cal 
O}_\Delta,{\cal O}_\Delta).
\end{eqnarray*}
\bigskip
\noindent{{\bf{\underbar{Proof:}}}} \ \ \ \ \   Let
\begin{equation}
R: \; K^._\alpha \to {\cal O}_{\Delta \cap U_\alpha} 
\tag{2.8}
\end{equation}

\noindent be the quasi isomorphism $R: \;K^0_\alpha \to {\cal O}_{\Delta \cap 
U_\alpha}$ the projection, and $R=0$
acting on $K^q_\alpha$, $q \not=0$, ${Hom}^.(K^.,K^.)$ is the total complex of 
the bicomplex $K^. \otimes
{\check K}^.$ whose first differential has cohomology concentrated in top 
degree, and so by $[TT2, (2.9)]$
$$
H^.({Hom}^.(K^.,K^.)) \xrightarrow[\cong]{R} H^.({Hom}^.(K^., {\cal O}_\Delta)) 
\cong \underbar {\it {Ext}}^._{{\cal O}_{X \times X}}({\cal O}_\Delta,{\cal 
O}_\Delta)
$$

\noindent Next $R$ also induces a chain map of twisted complexes
\begin{equation}
C^._{a,a}({\cal W}, {Hom}^.(K^.,K^.))\xrightarrow{R}C^._{a,1}({\cal 
W},{Hom}^.(K^.,{\cal O}_\Delta)) 
\tag{2.9}
\end{equation}

\noindent $R$ preserves the ${\check C}\!ech$ filtrations and the corresponding 
spectral
sequences. At $E_2$ level $R$ is just the identity map on $H^.(X\times X, 
\underbar
{\it {Ext}}^._{{\cal O}_{X \times X}}({\cal O}_\Delta,{\cal O}_\Delta))$ by 
preceding discussions.
Hence we have again by $[TT2, (2.9)]$
$$
H^.(C^._{a,a}({\cal 
W},{Hom}^.(K^.,K^.)))\xrightarrow[\cong]{R}H^.(C^._{a,1}({\cal
W},{Hom}^.(K^.,{\cal O}_\Delta)))\cong {Ext}^.(X \times X;{\cal O}_\Delta,{\cal 
O}_\Delta)
$$

\par The map $R$ in (2.8) can be given more detailed interpretations. By (2.6)

\begin{equation*}
R\;: \; {Hom}^q(K^., K^.) \to {Hom}(K^{-q}, {\cal O}_\Delta) \cong \bigwedge ^q 
T_X  
\end{equation*}

\noindent and by (2.1), $a^{i,-i+1}\bigm |_\Delta=0$ for $i \not= 1$. This shows 
that in fact we have a chain map of complexes
\begin{equation}
R\;:\;C^._{a,a}({\cal W},{Hom}^.(K^.,K^.)) \to \oplus C^.({\cal U}, 
\bigwedge^.T_X)
\tag{2.10}
\end{equation}

\noindent where the right hand side is the single complex associated to the 
${\check C}\!ech$
bicomplex with the differential in the coefficient sheaves $\wedge^.T_X$ being 
zero. Note
further that $R$ is obviously compatible with the cup products in the chain 
complexes. The
cup product in $C^._{a,a}({\cal W}, {Hom}^.(K^., K^.))$ is given by (1.10), the 
usual ${\check
C}\!ech$ cup product, together with compositions in coefficients. In cohomology 
this
corresponds to the Yoneda product in ${Ext}^.(X\times X; {\cal O}_\Delta,{\cal
O}_\Delta)(\rm{cf} [G,M])$. The restriction of this cup product to $\oplus 
C^.({\cal U},\wedge ^. T_X)$ via
$R$ clearly gives the ${\check C}\!ech$ cup product together with wedge product 
in the
$\wedge ^. T_X$.

\bigskip
\noindent{{\bf{\underbar{Theorem}} (2.11)}} \textit{The map (2.8) induces a 
chain map
preserving cup products
$$
R\;:\; C^._{a,a}({\cal W}, {Hom}^.(K^., K^.)) \to \oplus\,C^.({\cal U},\wedge 
^.T_X).
$$
\noindent In cohomology $R$ gives an isomorphism
$$
R\;:\;{Ext}^.(X\!\times\! X;{\cal O}_\Delta,{\cal O}_\Delta) \to 
\oplus\,H^.(X,\wedge ^.T_X)
$$
\noindent preserving the Yoneda product in ${Ext}^.$ and the cup product in 
$H^.(X,\wedge
^.T_X)$.}

\bigskip

\noindent{{\bf{\underbar{Proof:}}}} \ \ \ It only remains to point out that the 
same map $R$
is shown in (2.7) to be an isomorphism to its image.

\bigskip
\noindent{\bf \S3. Local Comparisons of Bar and Koszul Resolutions and Global 
consequences.}
\bigskip
\par Let $U \subset X$ be a Stein open set and $A=\Gamma(U,{\cal O})$ the ring 
of holomophic
functions on $U$. For any integer $q \ge 0$, let

\begin{equation}
{\cal B}^{-q}(A)=A \otimes_{\mathbb{C}}...\otimes_{\mathbb{C}} A=A^{\otimes 
(q+2)} 
\tag{3.1}
\end{equation}

\noindent ${\cal B}^{-q}(A)$ is an $A^e=A\otimes A$ module by multiplication in 
the first and last
factors. The bar resolution of $A$ is
$$
\to {\cal B}^{-q}(A)\xrightarrow{\partial}...\to {\cal 
B}^{-1}(A)\xrightarrow{\partial}{\cal B}^0(A) \to A \to 0
$$
\noindent where $\partial$ is the $A^e$ - linear map
$$
\partial(f_0\otimes...\otimes f_{q+1})=\sum^q_{i=0}(-1)^i f_0 \otimes...\otimes 
f_i
f_{i+1}\otimes...\otimes f_{q+1}
$$
\noindent On $U \times U$ we have the Koszul complex $K^._U$ defined in (1.1). 
We proceed to
construct a quasi isomorphism
$$
\Phi ^.\;:\; {\cal B}^. \to K^.
$$
\noindent Let $(z,\zeta)=(z^1,..,z^n,\zeta ^1,..,\zeta ^n)$ be coordinates on $U 
\times U$.
We will identify $A^e$ with $p^*_1 {\cal O}_U \otimes p^*_2 {\cal O}_U$.
$\Phi ^0\;:\; {\cal B}^0 \to K^0$ is the $A^e$ - linear map.
$$
\Phi ^0(f_0 \otimes f_1)= f_0(z)f_1(\zeta) \in K^0
$$
\noindent Next for $\eta=f_0 \otimes f_1\otimes f_2 \in {\cal B}^{-1}$, we set 
$$
\Phi^{-1}(\eta)=f_0 (z)\{P \Phi^0 \partial(1\otimes f_1\otimes 1)\}f_2(\zeta)
$$
\noindent where $P \;:\;K^0 \to K^{-1}$ is the ${\mathbb{C}}$ linear homotopy 
operator (1.2) on
$U\times U$, $\Phi^{-1}$ is $A^e$ linear by definition.

\par Assume inductively that $\Phi^{-k}$ is defined for $k \le {q-1}$, then set

\begin{equation}
\begin{aligned}
\Phi^{-q}(f_0 \otimes&...\otimes f_{q+1})\\
&=f_0 (z)\{ P \Phi^{-(q-1)}(\partial(1\otimes f_1 \otimes...\otimes f_q \otimes
1))\}f_{q+1}(\zeta)
\end{aligned}  
\tag{3.2}
\end{equation}
\noindent If $q>n$, set $\Phi^{-q}=0$. \newline
\bigskip

\noindent{\bf{\underbar{Lemma}(3.3)}}\ \ \ $\Phi^.\;:\; {\cal B}^. \to K^.$ 
\textit{is a quasi isomorphism}. 

\bigskip

\noindent{{\bf{\underbar{Proof:}}}} \ \ \ Both ${\cal B}^.$ and $K^.$ give 
resolutions of ${\cal
O}_{\Delta \cap U\times U}$ and clearly $\Phi^0$ induces identity map on ${\cal
O}_{\Delta \cap U\times U}$. To see that $\Phi^.$ commutes with differentials, 
let $\eta=f_0\otimes f_1 \otimes f_2 \in {\cal B}^{-1}$ \ \ \ then
$$
\Phi^0(\partial \eta)= f_0(z)(f_1(z)-f_1(\zeta))f_2(\zeta)
$$
\noindent On the other hand
\begin{equation*}
\begin{aligned}
d_K \Phi^{-1}(\eta)&=d_Kf_0(z)[P\Phi^0(\partial(1\otimes f_1\otimes 
1))]f_2(\zeta)\\
& =f_0(z)f_2(\zeta)(d_K P)(f_1(z)-f_1(\zeta)) \\
& =f_0(z)f_2(\zeta)(f_1(z)-f_1(\zeta))=\Phi ^0(\partial \eta).
\end{aligned}
\end{equation*}
\noindent where we have used the identity (1.3) on $K^0$:
$$
d_K\,P = 1- res
$$
\noindent and the fact $res(f_1(z)-f_1(\zeta))=0$.
\noindent Now assume $d_K \Phi^{-k}=\Phi ^{-(k-1)} \partial$ for $1 \le k\le 
{q-1}$.
\noindent Consider

\begin{equation}
\begin{aligned}
d_K \Phi ^{-q}&(f_0\otimes...\otimes f_{q+1})\\
& =f_0(z)f_{q+1}(\zeta)[d_K P \Phi ^{-(q-1)}(\partial(1\otimes f_1...\otimes 
f_g\otimes 1))] \\
& =f_0(z)f_{q+1}(\zeta)[(1-P d_K)\Phi^{-(q-1)}(\partial(1\otimes f_1...\otimes 
f_q \otimes 1))]
\end{aligned}
\tag{3.4}
\end{equation}

\noindent By inductive hypothesis $d_K \Phi ^{-(q-1)}=\Phi ^{-(q-1)}\partial$. 
So (3.4) reduces to
\begin{equation*}
\begin{aligned}
f_0(z)f_{q+1}(\zeta)\Phi ^{-(q-1)}&(\partial(1\otimes f_1...\otimes f_q \otimes 
1))\\
&=\Phi^{-(q-1)}(\partial(f_0\otimes ...\otimes f_{q+1})).
\end{aligned}
\end{equation*}
\noindent This finishes the proof of lemma.

\bigskip
\par To compare with the standard formula of the Hochschild-Kostant-Rosenberg 
isomophism, we
consider the map
\begin{equation}
\tilde \Phi^.\;:\;{\cal B}^. \xrightarrow{\Phi ^.} K^. \to K^.\bigm|_\Delta 
\tag{3.5}
\end{equation}
\noindent By (2.5), $K^.\bigm|_\Delta$ is naturally identified with $\Omega 
^._U$.

\bigskip
\noindent{\bf{\underbar{Lemma}(3.6)}} \ \ \ $\tilde \Phi^. (f_0 
\otimes...\otimes f_{q+1})= {1 \over
q!}f_0(z)f_{q+1}(z)df_1 \wedge...\wedge df_q$ \newline
\bigskip

\noindent{{\bf{\underbar{Proof:}}} \ \ \ We prove by induction on $q$. \newline
\noindent When $q=1$
$$
\Phi^{-1}(f_0 \otimes f_1\otimes f_2)=f_0(z)f_2(\zeta)P(f_1(z)-f_1(\zeta)).
$$
\noindent Hence the restriction to $\Delta$ gives
$$
\tilde \Phi^{-1}(f_0\otimes f_1 \otimes f_2)=f_0(z)f_2(z)df_1.
$$
\noindent Assume the lemma is valid for $k<q$. Then by (3.2):
\begin{equation}
\begin{aligned}
\Phi^{-q}(f_0\otimes &...\otimes f_{q+1})\\
&=f_0(z)f_{q+1}(\zeta) P\{\Phi ^{-(q-1)}[f_1 \otimes f_2...\otimes f_q \otimes 1 
-1 \otimes
f_1f_2 \otimes...\otimes 1 \\
&...\pm 1\otimes f_1..\otimes f_q]\} 
\end{aligned}
\tag{3.7}
\end{equation}

\noindent In (3.7) we expand again using (3.2) for $\Phi^{-(q-1)}$,
\begin{equation}
\begin{aligned}
f_0(z)f_{q+1}(\zeta) P & \{ f_1(z)P \Phi ^{-(q-2)} 
\partial(1\otimes f_2 \otimes ...\otimes f_q \otimes 1)- P \Phi^{-(q-2)} 
\partial (1\otimes f_1f_2\otimes...\otimes 1)\\
 & ....\pm f_q (\zeta)P\Phi^{-(q-2)}\partial(1\otimes f_1...\otimes 1) \}  
\end{aligned}
\tag{3.8}
\end{equation}

\noindent Using the fact that $P^2=0((1,3))$, and $P$ differentiates in $z$ 
variables and is therefore linear in functions in $\zeta$, it is clear that all 
the terms in the sum in (3.8) drop out except the first one. Upon restriction to 
$\Delta$, by induction hypothesis
$$
P\Phi^{-(q-2)} \partial(1\otimes f_2\otimes...\otimes f_q\otimes 1)\bigm 
|_\Delta ={1\over (q-1)!}\,df_2\wedge...\wedge \, df_q.
$$
\noindent Substitute this into (3.8) we get by (1.2)
$$
\tilde \Phi^{-q}(f_0\otimes...\otimes f_{q+1})= f_0(z)f_{q+1}(z){1\over {q!}}\, 
df_1 \wedge df_2\wedge...\wedge df_q.
$$
\noindent This proves the lemma.

\bigskip
\par Consider the bar complex ${\cal B}^.(A) \otimes_{A^e} A$ where $f_0\otimes 
f_1\otimes...f_q\otimes
f_{q+1}$ is identified with $f_0f_{q+1}\otimes f_1\otimes...\otimes f_q\otimes 
1$ in (3.1). We will denote
the bar complex by ${\cal B}^. \bigm |_\Delta$ or ${\cal B}^. \otimes {\cal 
O}_\Delta$. We also denote by $\tilde \Phi^.$ the restriction of the $A^e$ 
linear map (3.6) to ${\cal B}^. \bigm |_\Delta$.

\begin{equation}
\tilde \Phi^. \;:\; {\cal B}^. \bigm |_\Delta \to K^.\bigm|_\Delta 
\tag{3.9}
\end{equation}

\bigskip
\noindent The definition of ${\cal B}^.$ does not extend globally to make it 
a complex of sheaves on $X\times X$. 
However ${\cal B}^.\otimes {\cal O}_\Delta$ is well defined to 
give a global complex of sheaves on $X$, which is the 
standard Hochschild chain complex on $X$. 
By (3.6) the map (3.9) is just the map of 
Hochschild-Kostant-Rosenberg. 
${\cal B}^. \otimes {\cal O}_\Delta$ 
has a natural shuffle product $[W2, 9.4]$. 
For elements $f_0\otimes f_1\otimes...\otimes f_p$ 
and $f'_0 \otimes f_{p+1}\otimes...\otimes f_{p+q}$ their 
shuffle product $f_0\otimes f_1\otimes...\otimes f_p \# f'_0\otimes 
f_{p+1}\otimes...\otimes f_{p+q}$ is given by
\begin{equation}
\sum_{\scriptstyle\sigma \in \sum_{p+q} \atop\scriptstyle \sigma \,is \,a 
\,(p,q) \,shuf\!fle}\!(-1)^\sigma f_0 f'_0 \otimes 
f_{\sigma(1)}\otimes...\otimes f_{\sigma(p)}\otimes f_{\sigma(p+1)} 
\otimes...\otimes f_{\sigma (p+q)} 
\tag{3.10}
\end{equation}
\bigskip

\noindent It follows from $[W2,Cor 9.4.4]$ that (3.9) is compatible with the 
shuffle product in ${\cal B}^. \bigm |_\Delta$ and exterior product of forms in 
$K^.\bigm |_\Delta$. 

\bigskip
\par Combining the chain map $\Phi^. \;:\; {\cal B}^. \to K^.$ with the 
projection \newline
$R\;:\; K^. \to {\cal O}_\Delta$, (2.8), in second factor of ${Hom}^.(K^., 
K^.)$, we have a chain map

\begin{equation}
\Psi^.\;:\; {Hom}^.(K^.,K^.) \to {Hom}^.({\cal B}^., {\cal O}_\Delta)  
\tag{3.11}
\end{equation}

\noindent which in cohomology gives identity on the sheaves
$\underbar {\it {Ext}}^._{{\cal O}_{X \times X}}({\cal O}_\Delta,{\cal 
O}_\Delta)$. ${Hom}^.({\cal B}^., {\cal O}_\Delta)$ is the standard Hochschild 
cochain complex of $X$. Let $HH^.(X)$ be the hyperchomology of ${Hom}^. ({\cal 
B}^., {\cal O}_\Delta)[S,W1,WG]$. Using a Stein covering of $X$ the Hochschild 
cohomology of $X$ is given by the total cohomology of the ${\check C}\!ech - 
{Hom}^.$ bicomplex:

\begin{equation}
HH^.(X) = H^.(C^.({\cal U}, {Hom}^.({\cal B}^., {\cal O}_\Delta))
\tag{3.12}
\end{equation}

\bigskip
\par In the setting of the cochain complex ${Hom}^.({\cal B}^., {\cal 
O}_\Delta)$ the Hochschild-Kostant-
Rosenberg map (abbreviated HKR) is given by
\begin{equation}
\bigwedge ^. T_X \xrightarrow{HKR} {Hom}^.({\cal B}^., {\cal O}_\Delta). 
\tag{3.13}
\end{equation}

\bigskip
\noindent where $\upsilon_1\wedge...\wedge \upsilon_q \in \bigwedge ^q T_X$ acts 
on $f_0\otimes f_1\otimes...\otimes f_q$ by
$$
{1\over q!}\sum_{\sigma \in \sum_q} (-1)^\sigma f_0 \, v_{\sigma 
(1)}(f_1)...v_{\sigma (q)}(f_q).
$$

\bigskip
\noindent This is the dual of (3.9), therefore we get a commutative diagram:

\begin{equation}
\begin{array}{ccc}
{{Hom}^.(K^.,K^.)} & \xrightarrow{\displaystyle \Psi^.} & {{Hom}^.({\cal B}^., 
{\cal O}_\Delta)}\\
\searrow{R} & & \nearrow{HKR}\\
& {\wedge^. T_X}
\end{array}
\tag{3.14}
\end{equation}

\bigskip
\noindent By the remark following (3.10), HKR preserves products, and by \S2 $R$ 
preserves products, therefore $\Psi^.$ also preserves products.
\bigskip
\par Consider the composition on global ${\check C}\!ech$ complexes:
\begin{equation}
C^._{a,a}({\cal W}, {Hom}^.(K^.,K^.))\xrightarrow{R}C^._{a,1}({\cal W}, 
{Hom}^.(K^.,{\cal O}_\Delta))\xrightarrow{\displaystyle \Psi^.}C^.({\cal 
U},{Hom}^.({\cal B}^., {\cal O}_\Delta)) 
\tag{3.15}
\end{equation}

\bigskip
\noindent which we still denote by $\Psi^.$. Recall that $a^{i,-i+1} 
\bigm|_\Delta=0$ for $i\ge 2$. This implies that the composition $\Psi^.$ is a 
chain map. $\Psi^.$ induces isomorphism in total cohomology since it presserves 
spectral sequences induced by ${\check C}\!ech$ filtrations and it is 
isomorphism on $E_2$. We summarize the results in the following theorem.
\bigskip

\noindent{{\bf{\underbar{Theorem}} (3.16)}} \textit{There is a commutative 
diagram of chain complexes presserving cup products:}

\begin{equation*}
\begin{array}{ccc}
{C^._{a,a}({\cal W}, {Hom}^.(K^.,K^.))} & \xrightarrow{\displaystyle \Psi^.} & 
{C^.({\cal U},{Hom}^.({\cal B}^., {\cal O}_\Delta))}\\
\searrow{R} & & \nearrow{HKR}\\
& {C^.({\cal U}, \wedge^.T_X)}
\end{array}
\end{equation*}

\bigskip
\noindent \textit{The induced maps on cohomology are isomorphisms preserving cup 
products:}
\begin{equation*}
\begin{array}{ccc}
{{Ext}^.(X\!\times\!X;{\cal O}_\Delta,{\cal O}_\Delta)} & 
\xrightarrow{\displaystyle \Psi^.} & {HH^.(X)}\\
\searrow{R} & & \nearrow{HKR}\\
& {\oplus H^.(X,\wedge^.T_X)}
\end{array}
\end{equation*}

\bigskip

\noindent{\bf \S4. Hochschild Homology and $Tor$ Functors}
\bigskip
\par One way to represent the global $Tor$ functor $Tor.(X;{\cal F},{\cal G})$ 
for a pair of coherent sheaves ${\cal F}$ and ${\cal G}$ on $X$ is to construct twisted resolutions of ${\cal F}$ and ${\cal G}$ over a covering ${\cal 
U}$ of $X$, and then make use of the tensor product of twisting cochains $[OTT, 
\S 5]$. This involves more complicated constructions. For our present purpose 
which is to represent $Tor.(X\times X; {\cal O}_\Delta,{\cal O}_\Delta)$ and 
have an action by ${Ext}^.(X\times X;{\cal O}_\Delta,{\cal O}_\Delta)$, a 
simpler approach is possible.
\bigskip

\par Let $K^._\alpha$ be the Koszul complex from \S1. We set

\begin{equation}
\begin{aligned}
(K^._\alpha \otimes K^._\beta)^{-q}=\sum_{i+j=q} K^{-i}_\alpha \otimes 
K^{-j}_\beta 
&\cong \sum_{i+j=q} Hom({\check K}^j_\beta, K^{-i}_\alpha)={Hom}^{-q}({\check 
K}^._\beta, K^._\alpha)
\end{aligned} 
\tag{4.1}
\end{equation}

\bigskip
\noindent This tensor product has differential
$$
d_\otimes(\eta_\alpha \otimes \xi _\beta)= (d_K)_\alpha \eta_\alpha \otimes 
\xi_\beta + (-1)^{deg \, \eta}\eta_\alpha \otimes (d_K)_\beta \xi_\beta$$
$$=(d_K)_\alpha (\eta_\alpha \otimes \xi_\beta) + (-1)^{deg\, \eta+deg \, \xi + 
1}(\eta_\alpha \otimes \xi_\beta)(d_{\check K})_\beta
$$
\bigskip
\noindent Let $R\;:\;(K^. \otimes K^.) \to (K^. \otimes {\cal O}_\Delta)$ be the 
map induced by (2.8) in the second factor. As $d_K$ has homology concentrated at 
the top degree we have, by (2.5), over $U_\alpha \times U_\alpha \cap U_\beta 
\times U_\beta$

\begin{equation}
\underbar {\it {Tor}}_i^{{\cal O}_{X\times X}}({\cal O}_\Delta,{\cal O}_\Delta) 
\cong H^i((K^._\alpha \otimes K^._\beta)^.) \stackrel {R}{\cong}\Omega^i_\Delta 
\tag{4.2}
\end{equation}

\bigskip
\noindent Next to put a twisted differential on 
$$
C^.({\cal W},(K^.\otimes K^.)^.)\cong C^.({\cal W}, {Hom}^.({\check K}^., K^.))
$$

\noindent where ${\cal W}$ is a covering as in \S1, we will constuct a twisted 
resolution for ${\check K}$ and define $D_{a,{\check a}}$ as in $[OTT,\S1]$. A 
chain homotopy operator ${\check P}_\alpha$ on ${\check K}^._\alpha$ over the 
open set $U_\alpha \times U_\alpha$ is given by 
$$
{\check P}_\alpha (f{\check e}^I)=\sum^\eta_{j=1}\left(\int ^1_0 
t^{n-|I|}{\partial f \over \partial z^j_\alpha}(\zeta_\alpha +t(z_\alpha 
-\zeta_\alpha),\zeta_\alpha)dt \right)\iota_{e_j}{\check e}^I
$$
\noindent where $\iota$ denotes contraction.
\bigskip
\noindent It satisfies $([TT1, (9.19)])$

\begin{equation}
(d_{\check K})_\alpha {\check P}_\alpha + {\check P}_\alpha (d_{\check 
K})_\alpha= 1 - {\check {res}} 
\tag{4.3}
\end{equation}

$$({\check P}_\alpha)^2 =0$$
\bigskip
\noindent where ${\check {res}}(f{\check e}^1_1 \wedge...\wedge {\check 
e}^n)(z,\zeta)=f(z,z){\check e}^1\wedge...\wedge {\check e}^n$, ${\check 
{res}}|{\check K}^q=0$ if $q \not=n$. Since ${\check K}^.\bigm |_\Delta \cong 
\bigwedge^.T_X$ it follows that locally we have a resolution
$$
0 \to {\check K}^. \to \bigwedge ^nT_X \to 0
$$
\bigskip
\noindent Now using ${\check P}$ we construct a twisting cochain ${\check a}$ 
for ${\check K}^.$ by exactly the same recipe as the twisting cochain $a$ for 
$K^.(cf.[TT2,\S2])$. It should be noted that
$$
{\check a}^{i,-i+1}_{\alpha_0...\alpha_p} \in {Hom}^.({\check 
K}^._{\alpha_p},{\check K}^._{\alpha_0}) \cong \oplus {\check K}^._{\alpha_0} 
\otimes K^._{\alpha_p}
$$
\bigskip
\noindent and this is different from the one obtained by dualizing $a$ in 
$[TT2,(1.5),(1.7)]$, where the cochain takes value in 
$K^._{\alpha_0} \otimes {\check K}^._{\alpha_p}$. 
Here $C^r({\cal W},(K^. \otimes K^.)^{-s})$ are the ${\check C}\!ech\; 
r-cochains$
$$
c_{\alpha _0...\alpha _r} \in (K^._{\alpha _0} \otimes K^._{\alpha _r})^{-s}
$$
\bigskip
\noindent There is an action

\begin{equation}
\begin{aligned}
C^p({\cal W},{Hom}^q(K^.,K^.)) \times C^r({\cal W},(K^.\otimes K^.)^{-s})
&\to C^{p+r}({\cal W},(K^.\otimes K^.)^{-s+q}) 
\end{aligned}
\tag{4.4}
\end{equation}

\begin{equation*}
(f^{pq}\cdot c^{r,-s})_{\alpha _0...\alpha _{p+r}}=(-1)^{qr}f^{pq}_{\alpha 
_0...\alpha_p}c^{r,-s}_{\alpha_p...\alpha_{p+r}}
\end{equation*}

\noindent The differential $D_{a,{\check a}}$ on $C^.({\cal W},{Hom}^.({\check 
K}^.,K^.))$ is given by
$$
D_{a,{\check a}}c=\delta c+ a \cdot c+(-1)^{deg \,c+1}c \cdot {\check a}
$$
\bigskip
\noindent where $deg\, c$ is total degree as in (1.17). \newline
For $f\in C^.({\cal W},{Hom}^.(K^.,K^.))$ we have
$$
D_{a,{\check a}}(f\cdot c)= D_{a,a}(f) \cdot c+(-1)^{deg\, f}f \cdot 
D_{a,{\check a}}(c)
$$
\bigskip
\noindent $C^. _{a,{\check a}}({\cal W},{Hom}^. ({\check K}^.,K^.))$ denotes the 
singly graded complex with the differential $D_{a,{\check a}}$. \\
\noindent The global functor ${Tor}_m(X \times X;{\cal O}_\Delta,{\cal 
O}_\Delta)$ is represented by 
$$
{Tor}_m(X\! \times\! X;{\cal O}_\Delta,{\cal O}_\Delta)\cong 
H^{-m}\left(C^._{a,{\check a}}({\cal W},{Hom}^.({\check K}^.,K^.))\right)
$$
\bigskip
\noindent and arguments paralled to \S2 show that
\begin{equation}
\begin{aligned}
R\;:\;C^._{a,{\check a}}({\cal W},{Hom}^.({\check K}^.,K^.))\to C^._{1,{\check 
a}}({\cal W},{Hom}^.({\check K}^.,{\cal O}_\Delta))
&\cong \oplus C^.({\cal U},\Omega^._X) 
\end{aligned}
\tag{4.5}
\end{equation}

\noindent induces an isomorphism.
\begin{equation}
H^{-m}\left(C^._{a,{\check a}}({\cal W},{Hom}^. ({\check 
K}^.,K^.))\right)\xrightarrow{R} \bigoplus_i H^{i-m}(C^.({\cal U},\Omega^i_X)) 
\tag{4.6}
\end{equation}

\bigskip
\noindent which gives degeneration of spectral sequence for $Tor$:
$$
Tor_m(X\!\times\!X;{\cal O}_\Delta,{\cal O}_\Delta)\cong \bigoplus_i 
H^{i-m}(X\!\times\!X;\underbar{\it {Tor}}^{{\cal O}_{X\!\times\!X}}({\cal 
O}_\Delta,{\cal O}_\Delta))
$$
\bigskip
\par The action ${Hom}^q(K^.,K^.) \times {Hom}^{-s}({\check K}^.,K^.)\to 
{Hom}^{-s+q}({\check K}^.,K^.)$
\bigskip
\noindent gives in cohomology:
$$
\underbar {\it {Ext}}^q_{{\cal O}_{X\!\times\!X}}({\cal O}_\Delta,{\cal 
O}_\Delta)) \times \underbar{\it {Tor}}^{{\cal O}_{X\!\times\!X}}_s({\cal 
O}_\Delta,{\cal O}_\Delta)\to \underbar{\it {Tor}}^{{\cal 
O}_{X\!\times\!X}}_{s-q}({\cal O}_\Delta,{\cal O}_\Delta)
$$
\bigskip
\noindent or the contractions:
$$
\bigwedge^q\,T_X \times \Omega^s_X \to \Omega^{s-q}_X
$$
\bigskip
\noindent Globally the action (4.4) induces in cohomology
$$
{Ext}^k(X\!\times\!X;{\cal O}_\Delta,{\cal O}_\Delta) \times 
{Tor}_m(X\!\times\!X;{\cal O}_\Delta,{\cal O}_\Delta)$$
$$\to {Tor}_{m-k}(X\!\times\!X;{\cal O}_\Delta,{\cal O}_\Delta).
$$
\bigskip
\noindent Finally the Hochschild homology of $HH_m(X)$ is given by the total 
cohomology $[W1,WG]$:
\begin{equation}
HH_m(X)= H^{-m}(C^.({\cal U},{\cal B}^. \otimes {\cal O}_\Delta)) 
\tag{4.7}
\end{equation}

\bigskip
\noindent and via the chain map $\Phi ^.\;:\; {\cal B}^. \to K^.$ \ \ (\S2), 
 we have the isomorphism
 
\begin{eqnarray*}
H^{-m}(C^.({\cal U},{\cal B}^. \otimes {\cal O}_\Delta))&\xrightarrow{\sim}& 
H^{-m}(C^._{1,{\check a}}({\cal W},{Hom}^.({\check K}^.,{\cal O}_\Delta)))\\
&\xrightarrow{\simeq}& \bigoplus_i H^{i-m}(C^.({\cal U},\Omega^i_X)).
\end{eqnarray*}

\bigskip
\noindent By (3.6) this is just the standard HKR isomorphism on Hochschild 
homology. This also coincides with the map (4.6). 

\bigskip
\noindent We summarize in the following theorem.\\

\noindent{{\bf{\underbar{Theorem}} (4.8)}} 
\textit{The maps $\Phi$,(3.2), and $R$,(2.8), give a commutative diagram of 
chain maps.}
\bigskip
\begin{equation*}
\begin{array}{ccc}
{C^.({\cal U}, {\cal B}^.\otimes {\cal O}_\Delta)} & \xrightarrow{\displaystyle 
{\tilde \Psi^.}} & {C^._{1,{\check a}}({\cal W},{Hom}^.({\check K}^., {\cal 
O}_\Delta))}\\
\searrow{HKR} & & \swarrow{R}\\
& {C^.({\cal U}, \Omega^._X)} &
\end{array}
\end{equation*}

\noindent \textit{which induce isomophisms in cohomology:}

\begin{equation*}
\begin{array}{ccc}
{HH_m(X)} & \xrightarrow{\displaystyle {\tilde \Psi^.}} & 
{{Tor}_m(X\!\times\!X;{\cal O}_\Delta, {\cal O}_\Delta))}\\
\; \searrow{HKR}& & \swarrow{R}\\
& {\bigoplus_i H^{i-m}(X, \Omega ^i_X)} &
\end{array}
\end{equation*}

\textit{Furthermore the action (4.4)induces in cohomology}
$$
{Ext}^l(X\!\times\!X;{\cal O}_\Delta,{\cal O}_\Delta) 
\times{Tor}_m(X\!\times\!X;{\cal O}_\Delta,{\cal O}_\Delta)\to 
{Tor}_{m-l}(X\!\times\!X;{\cal O}_\Delta,{\cal O}_\Delta)
$$

\textit{which corresponds, via $R$, to the contractions.}
$$
H^j(X,\wedge^k T_X) \times H^p(X,\Omega^q_X) \to H^{j+p}(X, \Omega^{q-k}_X).
$$

\bigskip

\bibliographystyle{amsalpha}

Yue Lin L. Tong\\
Department of Mathematics\\
Purdue University\\
West Lafayette, IN. 47907\ \ USA\\

\smallskip
\noindent I-Hsun Tsai\\
Department of Mathematics\\
National Taiwan University\\
Taipei, Taiwan

\end{document}